\documentclass[11pt]{article}%
\usepackage{amssymb}
\usepackage{amsmath}
\usepackage{mitpress}
\usepackage{graphicx}
\usepackage{amsfonts}%
\providecommand{\U}[1]{\protect\rule{.1in}{.1in}}
\providecommand{\U}[1]{\protect\rule{.1in}{.1in}}
\providecommand{\U}[1]{\protect\rule{.1in}{.1in}}
\providecommand{\U}[1]{\protect\rule{.1in}{.1in}}
\newtheorem{theorem}{Theorem}[subsection]

\newtheorem{corollary}[theorem]{Corollary}

\newtheorem{lemma}[theorem]{Lemma}

\newtheorem{proposition}[theorem]{Proposition}
\newtheorem{remark}[theorem]{Remark}

\newenvironment{proof}[1][Proof]{\textbf{#1.} }{\ \rule{0.5em}{0.5em}}
\newdimen\dummy
\dummy=\oddsidemargin
\ifx\pdfoutput\relax\let\pdfoutput=\undefined\fi
\newcount\msipdfoutput
\ifx\pdfoutput\undefined\else
\ifcase\pdfoutput\else
\ifx\paperwidth\undefined\else
\ifdim\paperheight=0pt\relax\else\pdfpageheight\paperheight\fi
\ifdim\paperwidth=0pt\relax\else\pdfpagewidth\paperwidth\fi
\fi\fi\fi
\ifx\pdfoutput\relax\let\pdfoutput=\undefined\fi
\newcount\msipdfoutput
\ifx\pdfoutput\undefined\else
\ifcase\pdfoutput\else
\ifx\paperwidth\undefined\else
\ifdim\paperheight=0pt\relax\else\pdfpageheight\paperheight\fi
\ifdim\paperwidth=0pt\relax\else\pdfpagewidth\paperwidth\fi
\fi\fi\fi
\begin{document}

\title{\textsc{Special biserial coalgebras and representations of quantum SL(2)}}
\author{William Chin\\DePaul University, Chicago Illinois 60614}
\maketitle

\begin{abstract}
We develop the theory of special biserial and string coalgebras and other
concepts from the representation theory of quivers. These tools are then used
to describe the finite dimensional comodules and Auslander-Reiten quiver for
the coordinate Hopf algebra of quantum SL(2) at a root of unity. We also
compute quantum dimensions and the stable Green ring.

\end{abstract}

Let $C=\Bbbk_{\zeta}[SL(2)]$ denote the quantized coordinate Hopf algebra of
$SL(2)$ as defined in [APW] and [Lu2] at a root of unity of odd order over a
field of characteristic zero. In this article we study the category of
finite-dimensional comodules, a category that is equivalent to the category of
finite-dimensional modules over (a suitable quotient of) the quantized
hyperalgebra $U_{\zeta}$. Our approach uses representations of the Gabriel
quiver associated to $C$, methods from the representation theory of quivers,
and most notably, the theory of string algebras and special biserial algebras.
The methods of this paper demonstrate the utility of coalgebraic methods,
which we expect will be applicable to other comodule categories.

We discuss some required coalgebra representation theory in section 1 and then
use the more general results to concentrate on the case of the quantized
coordinate Hopf algebra $C$ at a root of unity over a field of characteristic
zero in section 2. In [CK] we determined the structure of the injective
indecomposable comodules along with the block decomposition of $C$. Here we
completely determine the finite-dimensional indecomposable comodules, the
Auslander-Reiten quiver and almost split sequences.

We introduce the notions of a special biserial coalgebra and string coalgebra,
the latter modifying the definition in [Sim2], [Sim3]. These notions are
coalgebraic versions of machinery known for algebras [Ri], [BR], [Erd]. It
turns out that $C$ is special biserial and that its representations are
closely related to an associated string subcoalgebra, and this allows a
listing of indecomposables and a description of the Auslander-Reiten quiver.
One consequence is that the coalgebra $C$ is of tame of discrete comodule type
[Sim2] and is the directed union of finite-dimensional coalgebras of finite
type. This is in contrast with the situation for restricted representations,
see [Xia], [Sut], [RT]. We compute the Green ring modulo
projective-injectives, or \textit{stable Green ring}, in which the syzygy
operation is given by multiplication by the isomorphism class of a certain
Weyl module. Here again the structure is simpler than for restricted
representations (cf. [EGST]); the stable Green ring is a polynomial ring
$\mathbb{Z}[x,y,w^{\pm1}]$ in three indeterminates modulo a rescaled Chebyshev
polynomial of the second kind in $x$.

The notion of quantum dimension for comodules over $C\ $was studied in [And].
Quantized tilting modules were of use in [Os], where module categories over
the tensor category of finite-dimensional comodules $\mathcal{M}_{f}^{^{C}}$
of $C$ were characterized. Here we compute the quantum dimensions for
finite-dimensional $C$-comodules using string modules and the Lusztig tensor
product theorem. As a result we find that the only comodules of quantum
dimension zero are the nonsimple injective comodules, a fact observed for
quantized tilting modules in [And].

Let us provide an outline of the contents of this article. The first section
contains relevant coalgebra theory. We begin in 1.1 by reviewing the notions
of path coalgebras, basic coalgebras and the dual Gabriel theorem. We
introduce special biserial and string coalgebras in 1.2 and their
representation theory in 1.3. Section 2 begins with a summary of facts
concerning the coordinate Hopf algebra $C=\Bbbk_{\zeta}[SL(2)]$ at a root of
unity where $\Bbbk$ is of characteristic zero. In 2.1 we compute the basic
coalgebra of $C$ as a subcoalgebra of the path coalgebra of its Gabriel
quiver. Using string coalgebra theory we list the finite-dimensional
indecomposable $C$-comodules in 2.2. In 2.3 we show that every simple comodule
is in the syzygy-orbit of a simple comodule by direct computation. This yields
the description of almost split sequences, and the AR quiver in 2.4. The
stable Green ring and quantum dimensions are computed in 2.5 and 2.6.

\bigskip

\ \textbf{Notation }Let $C$ denote a coalgebra over the fixed base field
$\Bbbk$. Set the following

\bigskip

$\mathcal{M}^{C}$ the category of right $C$-comodules

$\mathcal{M}_{f}^{^{C}}$ the category of finite-dimensional right $C$ -comodules

$\mathcal{M}_{q}^{^{C}}$ the category of quasifinite right $C$ -comodules

ind$(C)$ the category finite-dimensional indecomposable right $C$-comodules

\bigskip

$h_{-C}(\_,\_)$ the cohom functor

$\square$ the cotensor (over $C$)

\textrm{D} the $\Bbbk$-linear dual $\operatorname*{Hom}_{\Bbbk}(\_,)$

$(-)^{\ast}$ the functor $h_{-C}(\_,C):\mathcal{M}_{q}^{^{C}}\rightarrow
\mathcal{M}^{C}$

$A^{0}$ the finitary dual coalgebra of the algebra $A.$

\bigskip

The author wishes to acknowledge the advice of K. Erdmann concerning special
biserial algebras and D. Simson for his comments.

\section{Some coalgebra theory}

\subsection{Quivers and path coalgebras}

Green [Gr] (also see [Ch]) showed that the basic structure theory for
finite-dimensional algebras carries over to coalgebras, with injective
indecomposable comodules replacing projective indecomposables. Define the
(\textit{Gabriel- }or\textit{\ Ext-}) \textit{quiver} of a coalgebra $C$ to be
the directed graph $Q(C)$ with vertices $\mathcal{G}$ corresponding to
isoclasses of simple comodules and $\mathrm{dim}_{\Bbbk}\mathrm{Ext}^{1}(h,g)$
arrows from $h$ to $g$, for all $h,g\in\mathcal{G}$. The blocks of $C$ are
(the vertices of) components of the undirected version of the graph $Q(C)$. In
other words, the blocks are the equivalence classes of the equivalence
relation on $\mathcal{G}$ generated by arrows. The indecomposable or
\textquotedblleft block\textquotedblright\ coalgebras are the direct sums of
injective indecomposables having socles from a given block.

Write $C\in\mathcal{M}^{C}$ as a direct sum of indecomposable injectives with
multiplicities. Let $E$ denote the direct sum of the indecomposables where
each indecomposable occurs with multiplicity one. The comodule $E$ is called a
\textquotedblleft basic\textquotedblright\ injective for $C$. $E$ can be
written as $C\leftharpoonup e$ where $e\in DC$ is an idempotent in the dual
algebra. Define the \textit{basic coalgebra} to be the coendomorphism
coalgebra
\[
B=B_{C}=h_{-C}(E,E).
\]
Basic coalgebras were first constructed by Simson [Sim] and were rediscovered
in [CMo] and [Wo]. As observed in [CG], $B$ is more simply described as the
coalgebra $e\rightharpoonup C\leftharpoonup e$ (or just $eCe$), which is a
noncounital homomorphic image of $C$. Of course $B$ is Morita-Takeuchi
equivalent to $C$, meaning that their comodule categories are equivalent. If
$\Bbbk$ is algebraically closed, then $B$ is a pointed coalgebra.

Let $Q$ be a quiver (not necessarily finite) with vertex set $Q_{0}$ and arrow
set $Q_{1}.$ For a path $p$, let $s(p)$ denote the start (or source) of $p$
and let denote $t(p)$ its terminal (or target). The \textbf{path coalgebra}
$\Bbbk Q$ of $Q$ is defined to be the span of all paths in $Q$ with coalgebra
structure
\[%
\begin{array}
[c]{c}%
\Delta(p)=\sum_{p=p_{2}p_{1}}p_{2}\otimes p_{1}+t(p)\otimes p+p\otimes s(p)\\
\varepsilon(p)=\delta_{|p|,0}%
\end{array}
\]
where $p_{2}p_{1}$ is the concatenation $a_{t}a_{t-1}...a_{s+1}a_{s}%
...a_{1}\,$of the subpaths $p_{2}=a_{t}a_{t-1}...a_{s+1}$ and $p_{1}%
=a_{s}...a_{1}$ ($a_{i}\in Q_{0}$)$.$ Here $|p|=t$ denotes the length of $p$
and the starting vertex of $a_{i+1}\,$ is required to be the end of $a_{i}.$
The paths $p_{2}$ occurring in the sum are subpaths of\ $p$ with the same
terminal $t(p).$ The span of these subpaths (together with $t(p)$) span the
coideal (i.e. subcomodule) generated by $p.$ Similarly, any subcoalgebra
containing a path $p$ contains all subpaths.

Thus $vertices$ are group-like elements, and if $a$ is an arrow $g\leftarrow
h$, with $g,h\in Q_{0},$ then $a$ is a $(g,h)$- skew primitive, i.e., $\Delta
a=g\otimes a+a\otimes h.$ It is apparent that $\Bbbk Q$ is pointed with
coradical $(\Bbbk Q)_{0}=\Bbbk Q_{0}\,$and the degree one term of the
coradical filtration is $(\Bbbk Q)_{1}=\Bbbk Q_{0}\oplus\Bbbk Q_{1}.$

A subcoalgebra $B$ of the path coalgebra $\Bbbk Q$ is said to be
\textit{admissible} if $B$ contains $Q_{0}$ and $Q_{1}.$ In this case
$Q(B)=Q.$ There is a coalgebraic version (with no finiteness restrictions)\ of
a fundamental result of Gabriel for finite dimensional algebras:

\begin{theorem}
([CMo][Wo]) Every coalgebra $C$ over an algebraically closed field is
Morita-Takeuchi equivalent to an admissible subcoalgebra of $\Bbbk Q(C)$.
\end{theorem}

When $b=\sum a_{i}p_{i}$ $\in$ $B$ is a $\Bbbk$-linear combination of distinct
paths $p_{i}$ with $a_{i}\neq0$ for all $i$, we say that $b$ is a
\textit{reduced element}, and we denote the support of $b$ by
\[
\mathrm{supp}(b)=\{p_{i}\}.
\]
We shall also say that $p_{i}$ \textit{appears in} $B$ and write
\[
p_{i}\vdash B
\]
if $p_{i}$ $\in$\textrm{supp}$(b)$ for some $b\in B.$ Clearly $p\vdash B$ if
and only if $p\in B$ for all paths $p,$ exactly in the case where $B$ is the
span of a set of paths.

\subsection{Special biserial and string coalgebras}

\begin{description}
\item[Definition] A coalgebra $C$ is said to be a \textbf{special biserial
coalgebra} if its basic coalgebra is Morita-Takeuchi equivalent to an
admissible subcoalgebra $B\subset\Bbbk Q$ of its Gabriel quiver $Q$ such that
\end{description}

\begin{enumerate}
\item[S1] Every vertex of $Q$ is the start of at most two arrows and the end
of at most two arrows.

\item[S2] Given any arrow $b$, there is at most one arrow $a$ such that
$ab\vdash$ $B$ and at most one arrow $c$ such that $bc$ $\vdash$ $B.$
\end{enumerate}

If in addition,

\begin{enumerate}
\item[S3] $B$ is spanned by a set of paths (containing the vertices and arrows
of $Q$).
\end{enumerate}

we say that $C$ is a \textbf{string coalgebra.}

\bigskip

Note that the terms are defined up to Morita-Takeuchi equivalence and depend
on the embedding of the basic coalgebra into its path coalgebra. When we say
that $B\subset\Bbbk Q$ is a special biserial coalgebra, we adopt the
convention that $B$ is admissible and is special biserial with respect to
$\Bbbk Q.$ String algebras and special biserial algebras are defined dually,
given a possibly infinite bound quiver $(Q,I)$, see [BR][Erd][SW]. String
coalgebras have been defined for coalgebras specified by an ideal of relations
in the path algebra in [Sim3].

We make a few more elementary remarks. Let $B\subset\Bbbk Q$ be a coalgebra.
We let $I(g)$ be an injective indecomposable in $\mathcal{M}^{B}$, having
socle $S(g)$ corresponding to the vertex (group-like) $g\in Q_{0}.$ It is
useful to note that an injective envelope of the simple right $\Bbbk Q$
comodule $S(g)$ is $\Bbbk Q\leftharpoonup e_{g}$ where $e_{g}$ is the
idempotent dual to $g$ in the dual algebra $D(\Bbbk Q),$ using the right hit
action. It is easy to see that $\Bbbk Q\leftharpoonup e_{g}$ is simply the
span of all paths ending at $g$, and also that $I(g)=$ $(\Bbbk Q\leftharpoonup
e_{g})\cap B=Be_{g}.$

\bigskip

If $B$ is special biserial, then we can be specific about the form of $I(g)$.

\begin{proposition}
Let $B\subset\Bbbk Q$ be a special biserial coalgebra and let $I(g)$ be a
finite-dimensional injective indecomposable right $B$-comodule. Then one of
the following cases holds:\newline(a) $I(g)$ is a uniserial comodule,
generated by a path.\newline(b) $I(g)$ is the sum of two distinct uniserial
comodules, each generated by a path, where the only common vertex of the two
paths is $g.$ \newline(c) $I(g)$ is the sum of two distinct uniserial
comodules, each generated by a path, where the only common vertices of the two
paths is $g$ and a common starting vertex. \newline(d) $I(g)$ is generated by
$p+\lambda q,$ $0\neq\lambda\in\Bbbk\mathrm{,}$ where $p$ and $q$ are distinct
parallel paths, both ending at $g,$ with the same starting vertex, and no
other common vertices.
\end{proposition}

\begin{proof}
Suppose there is only one arrow only ending at $g$. \ Then condition S2 in the
definition implies (a).

Now assume that there are two arrows ending at $g$. If there are no common
vertices other than $g$, we are in case (b). Otherwise there are two paths
$p,q\vdash B$ in $I(g)$ having a common vertex other than $g$, and then
condition S2 forces it to be the common starting vertex, as in (c) and (d). If
$I(g)$ is spanned by paths, then $I(g)$ is as in (c).

Consider the case where $I(g)$ is not spanned by paths. In this case $B$
contains a linear combination, say $b=\sum\lambda_{i}p_{i}\in B,$ of at least
two paths not in $I(g),$ where each path involved has terminal vertex $g$.
Conditions S1 and S2 imply that the support $\{p_{i}\}$ of this element
consists of exactly two paths (say) $p$ and $q$ with $p\notin B$ and $q\notin
B$. Thus $I(g)$ contains an element $p+\lambda q$ as in the statement. It
follows easily that $s(p)=s(q)$ (otherwise we obtain the contradiction
$s(p)\rightharpoonup(p+\lambda q)=p\in B$), and as already stated, that $s(p)$
is the only common vertex of $p$ and $q.$ It remains to see that $I(g)$ is
generated by $p+\lambda q.$ If $b^{\prime}\in I(g)$ is a linear combination of
paths not in $B,$ with support size greater than $1$, then it follows as just
argued that the supp$(b)$ consists of exactly two paths. Furthermore these two
paths must be precisely $p$ and $q.$ This forces $b^{\prime}$ to be a scalar
multiple of $b.$ On the other hand, if $r\in I(g)$ is a path, then $t(r)=g$
and S1 forces the terminal arrow of $r$ to equal the terminal arrow of $p$ or
$q$ (unless $r=g$). Now S2 forces $r$ to be a proper subpath of $p$ or $q.$
This means that, say, $p=rs$ for a nonntrivial path $s$, and thus that
$s\rightharpoonup p=r\in I(g).$ This shows that $I(g)$ is generated by
$p+\lambda q,$ completing the proof of (d).
\end{proof}

\bigskip

It is easy to see that all cases in the Proposition above can occur. By
duality and the symmetry in the definition of special biserial coalgebras, it
follows that the comodules in case (d) of the Proposition above are
finite-dimensional and projective. In cases (b) and (c), the top of the
injective comodule is \textit{spanned} by the two generating paths (and thus
is not projective). In case (d), the top is spanned by the image of the
reduced element $p+\lambda q$. The only other case where $I(g)$ might be
projective is in the finite-dimensional uniserial case, where $I(g)$ is
generated as a comodule by a single path, as in part (a). It is a simple
matter to construct examples of such uniserial injectives.

We next point out that the requirement of admissibility in the definition of
special biserial and string coalgebra is superfluous.

\begin{lemma}
Suppose $B\subset\Bbbk Q$ is a subcoalgebra of a path coalgebra $\Bbbk Q$
satisfying S1 and S2. Then $B$ is a special biserial coalgebra.
\end{lemma}

\begin{proof}
Let $B$ be as in the statement. Write $B=B_{0}\oplus I,$ where $I=I_{1}\oplus
J$ is a coideal complement in $B$ for the coradical $B_{0}$, with $I_{1}=\Bbbk
Q_{1}\cap B$ and where $J$ is a $B_{0},B_{0}$-bicomodule complement for
$I_{1}$ in $I$. We will show that $B$ is embeddable as an admissible
subcoalgebra of $\Bbbk Q^{\prime}$ where $Q^{\prime}$ is the subquiver of $Q$
with vertices corresponding to the group-likes of $B$ and whose arrows are
selected as in the next paragraph. We will then argue that S1 and S2 hold for
the image of $B$ in $\Bbbk Q^{\prime}.$

Consider the skew-primitive spaces $I_{1}(x,y)=e_{x}I_{1}e_{y}$ for vertices
$x,y\in Q_{0}^{\prime}.$ Now according to Proposition 1.2.1, there are the two
nontrivial cases: (1) $I_{1}(x,y)$ is spanned by (one or two) arrows from $x$
to $y$, and (2) $I_{1}(x,y)$ is the span of a single linear combination of two
arrows, say $a+\lambda b$ , $0\neq\lambda\in\Bbbk,$ where $\{a,b\}=Q_{1}%
(x,y).$ Let $Q^{\prime}$ be the quiver whose arrows are the arrows in
$I_{1}(x,y)$ when $I_{1}(x,y)$ is spanned by arrows, and \textit{a choice of}
an arrow $a\vdash I_{1}(x,y)$ in case $I_{1}(x,y)$ is spanned by $a,b$ as in
case (2).

Now we embed $B$ into $\Bbbk Q^{\prime}$ as an admissible subcoalgebra by a
standard embedding of $B$ into the path coalgebra $\Bbbk Q$ as follows ([Ni],
also [CMo][Sim3][Wo]). We define the map $\pi_{1}:B\rightarrow\Bbbk
Q_{1}^{\prime}$ by projecting first onto $I_{1},$ and then as the identity on
components $I_{1}(x,y)$ in case (1) and sending a basis vector $a+\lambda b$
to the chosen arrow $a$ in case (2). Also we have the projection $\pi
_{0}:B\rightarrow B_{0}$ along $I$. Define a coalgebra embedding
$\theta:B\rightarrow\Bbbk Q^{\prime}$ by
\[
\theta(d)=\pi_{0}(d)+\sum_{n\geq1}\pi_{1}^{\otimes n}\Delta_{n-1}(d)
\]
for all $d\in B$. It is clear that $\theta$ embeds $B$ as an admissible
subcoalgebra $B^{\prime}=\theta(B)\subset$ $\Bbbk Q^{\prime}.$

It is not difficult to see that $B^{\prime}$ is a special biserial coalgebra.
Condition S1 is immediate from the hypothesis. For arrows appearing in case
(1), S2 is also immediately seen to be satisfied. We only need to check S2 for
the arrows $a$ appearing in case (2). In this case, by Proposition 1.2.1,
there can be \textit{no} paths of the form $ac$ or $ca$ for arrows $c$ such
that $ac\vdash B^{\prime}$ or $ca\vdash B^{\prime},$ as the arrow $a$
generates the injective-projective right submodule $\Bbbk x+\Bbbk a$ of $\Bbbk
Q^{\prime}.$ Thus S2 holds for all arrows in $Q^{\prime}.$ We have shown that
$B$ is isomorphic to an admissible subcoalgebra of $\Bbbk Q^{\prime}$
satisfying S1 and S2, completing the proof of the Lemma.
\end{proof}

For a comodule $M$, we let $\mathrm{rad}(M)$ denote the intersection of all
maximal subcomodules of $M$. Since submodules of rational $DC$-modules are
rational, this coincides with the standard usage. The next lemma may be well-known.

\begin{lemma}
Let $C$ be a coalgebra with finite-dimensional right C-comodule $M$. Let
$J=\mathrm{rad}(DC)\ $be the Jacobson radical of the algebra $DC.$ Then
$JM=\mathrm{rad}(M)$
\end{lemma}

\begin{proof}
It is elementary that $JM\subseteq\mathrm{rad}(M)$, so we shall prove the
other inclusion. Obviously $M/JM$ is annihilated by $J$. Therefore since
$M/JM$ is a rational $DC$-module, it follows that the coefficient space
cf$(M/JM)\subseteq J^{\perp}=\mathrm{corad}(C)$. Thus $M/JM$ is a completely
reducible $C$-comodule, so we conclude that $\mathrm{rad}(M)\subseteq JM,$ as
desired.\bigskip
\end{proof}

The following lemma is a coalgebraic version of results that seem to be
well-known for finite-dimensional algebras, cf. [Erd, II.1.3]. It is used in
the proof of the subsequent result.

\begin{lemma}
Let $C=P\oplus Q$ be a coalgebra where $P$ and $Q$ are injective right
subcomodules and $P$ is the direct sum of finite-dimensional indecomposable
projective right comodules.\newline(a) Assume $Q$ has no projective summands.
Then $\mathrm{rad}P\oplus Q$ is a subcoalgebra of $C$.\newline(b) Let $M$ be
an indecomposable right $C$-comodule, and assume $P$ is indecomposable with
$\mathrm{\operatorname*{cf}}(M)\nsubseteq\mathrm{rad}P\oplus Q.$ Then $M\cong
P.$
\end{lemma}

\begin{proof}
(a) Let $C^{\prime}=$rad$P\oplus Q$ and suppose that $C^{\prime}$ is not a
subcoalgebra of $C.$ Then $C^{\prime}$ is not a left $C$-comodule, so it not a
right $DC$-submodule and there exists $f\in DC$ such that $C^{\prime
}\leftharpoonup f\nsubseteq C^{\prime}.$ This right hit action extends to a
right comodule (and left $DC$-module) map $\phi:C\rightarrow C$ such that
$\mathrm{Im}(\phi)\nsubseteq C^{\prime}.$ Write $P=\bigoplus_{i}P_{i}$ and
$Q=\bigoplus_{j}Q_{j}$ for the respective Krull-Schmidt decompositions of $P$
and $Q$. By the previous lemma, we have $\phi($rad$P)=J\phi(C)\subset
\mathrm{rad}C\subset C^{\prime}.$ By projecting to $P,$ we see that there is a
comodule map from $Q$ to $P$ whose image is not contained in \ \textrm{rad}%
$P=\bigoplus_{i}$\textrm{rad}$P_{i}$. therefore, for some indices $i,j,$ there
exists a nonzero projection of $Q_{i}$ onto $P_{j}.$ This implies that $Q$ has
projective summand $P_{j}$. This concludes the proof of (a).

(b) We shall show that there exists $f\in\operatorname*{Hom}_{\mathrm{\Bbbk}%
}(M,\Bbbk)\cong\operatorname*{Hom}^{C}(M,C)$ such that $f$ induces a
surjection $M\rightarrow P.$ Indeed, the hypothesis yields $m\in M$ with
$\rho(m)=\sum m_{i}\otimes c_{i}$ with the $m_{i}$ linearly independent and
$c_{j}\notin$rad$P\oplus Q$ for some $j$. Let $f=m_{j}^{\ast}$, in a dual
basis for $\operatorname*{Hom}_{\mathrm{\Bbbk}}(M,\Bbbk)$ obtained by
extending the $m_{i}$ to a basis for $M.$ Then the right $C$-comodule mapping
\[
f^{\prime}=(1\otimes f)\rho:M\rightarrow C
\]
corresponds to $f$, and we have $f^{\prime}(m)=c_{j}.$ Composing with the
projection of $C$ onto $P$ yields a surjection of $M$ onto $P,$ since rad$P$
is the unique maximal subcomodule of $P$. Since $M$ is indecomposable, this
forces $M\cong P$ as desired. This completes the proof of the Lemma.
\end{proof}

The representation theory of special biserial coalgebras is closely related to
that of string coalgebras. If $C$ is special biserial, we can
\textquotedblleft remove\textquotedblright\ the reduced linear combinations of
paths and pass to a slightly smaller string coalgebra $C^{\prime},$ losing
only the obvious representations in the process. To be precise, let
$\mathcal{Q}=$ $\{g\in Q_{0}|I(g)$ is projective\} and put
\[
C^{\prime}=%
%TCIMACRO{\dbigoplus _{g\in\mathcal{Q}}}%
%BeginExpansion
{\displaystyle\bigoplus_{g\in\mathcal{Q}}}
%EndExpansion
\mathrm{rad}I(g)\oplus\,%
%TCIMACRO{\dbigoplus _{h\notin\mathcal{Q}}}%
%BeginExpansion
{\displaystyle\bigoplus_{h\notin\mathcal{Q}}}
%EndExpansion
I(h).
\]

\begin{proposition}
Let $B\subset\Bbbk Q$ be a special biserial coalgebra.\newline(a) Then
$B^{\prime}$ is a string coalgebra; in addition: $\mathrm{ind}$%
($B)=\mathrm{ind}(B^{\prime})\cup\{I(g)|$ $g\in\mathcal{Q}\}$ \newline(b) The
AR quiver of $B$ is obtained from the AR quiver for $B^{\prime}$ by attaching
an arrow from $\mathrm{rad}I(g)$ to $I(g)$ and an arrow from $I(g)$ to
$I(g)/S(g),$ for all $g\in\mathcal{Q}$.
\end{proposition}

\begin{proof}
Let $B$ be as in the hypothesis. Inspection of the list of possible forms for
the $I(g)$ in Proposition 1.2.1 shows that the ones of type (d) are exactly
the injective-projectives which are not uniserial and which are not spanned by
paths. It is easy to see that the radical of these $I(g)$'s are spanned by the
proper subpaths of $p$ and $q$ ending at $g$. The uniserial projectives in
case (a) which might occur are generated by a single path, and clearly their
radicals are also spanned by the proper subpaths. Thus we see that $B^{\prime
}$ is spanned by paths, and therefore that $B^{\prime}$ is a string coalgebra.

Part (b) of the preceding lemma shows that the $I(g),$ $g\in\mathcal{Q},$ are
the only indecomposables not already in $\mathrm{ind}(B^{\prime}).$ It remains
to prove that $B^{\prime}$ is a subcoalgebra of $B$. This assertion follows
directly from part (a) of the preceding lemma.

Part (b) follows from a standard argument (see e.g. [ARS], p. 169) that shows
that the only almost split sequences with a projective-injective occurring as
a summand of the middle term is of the form
\[
0\rightarrow\mathrm{rad}I(g)\rightarrow\mathrm{rad}I(g)/S(g)\oplus
I(g)\rightarrow I(g)/S(g)\rightarrow0
\]
This concludes the proof.
\end{proof}

\begin{remark}
One can replace
\[
\mathcal{Q}=\{g\in Q_{0}|I(g)\text{ is finite-dimensional projective}\}
\]
by
\[
\{g\in Q_{0}|I(g)\text{ is finite-dimensional projective and not uniserial}\}
\]
to obtain a larger string coalgebra in $B$. (cf. [Erd], II.1.3, II.6.4). If
$B$ is self-projective (i.e. quasi-coFrobenius, see [DNR]) to begin with, then
we obtain simply rad($B$) as a string coalgebra.
\end{remark}

\subsection{Representations of string coalgebras}

Consider a semiperfect string coalgebra $C$. By the theory of string algebras
mentioned in the introduction, the category of finite-dimensional comodules
$\emph{M}_{f}^{C}$ consists of \textit{string} modules and \textit{band}
modules which are locally nilpotent as quiver representations. Evidently, the
semiperfect assumption is not needed for the classification of
indecomposables, though it is relevant to the discussion of the AR quiver.

Here is a very basic example. Let $Q$ be a single loop with vertex $g$ and
arrow $\alpha$. Then $\Bbbk Q$ is a string coalgebra. However, the path
algebra $\Bbbk^{Q}=\Bbbk\lbrack x]$ is not a string algebra (as defined in
[BR], p. 157) because of the assumption that amounts to the coalgebra being
semiperfect. The finite dimensional modules are given by string and band
modules. $\ $The band modules corresponding to the band $\alpha$ and nonzero
scalar parameters $\lambda\in\Bbbk$ do not correspond to comodules. These are
the indecomposable modules annihilated by some power of $(x-\lambda),$
$0\neq\lambda\in\Bbbk$.

For the main example concerning quantum $SL(2)$ in this article in section 2
below, there are no bands, and we shall only be concerned with string modules,
which we presently describe. The string modules that will arise in section 2
will be relatively simple to describe, but for completeness we give the
general set-up.

\begin{itemize}
\item \textit{Letters} are arrows in $Q_{1}$ or formal inverses of arrows,
denoted by $\alpha^{-1}$ for $\alpha\in Q_{1}.$ We also set $(\alpha
^{-1})^{-1}=\alpha,$ $s(\alpha^{-1})=e(\alpha)$ and $e(\alpha^{-1})=s(\alpha
)$. A \textit{word} is defined to be a sequence $w=w_{n}...w_{1}$ of letters
with $e(w_{i})$ $=s(w_{i+1})$ and $w_{i}\neq w_{i+1}^{-1}$ for $0\leq i<n$,
$n\geq0$. Note that the empty word of length $n=0$ is allowed. The formally
inverted word $w^{-1}$ is defined to be $w_{1}^{-1}...w_{n}^{-1}$.

\item \textit{Strings} are defined to be the equivalence class of words where
each subpath (or its inverse) is in $B$, under the relation that identifies
each word with its inverse.

\item While we will not need them in section 2 (except to say that there are
not any of them), we mention that \textit{bands} are defined to be equivalence
classes of words, under cyclic permutation, whose powers $w^{m}$, $m\geq1$ are
defined (but are not themselves powers of words) and such that every subpath
of every power $w^{m}$ is in $B$. One-parameter families of modules are
associated to bands ([Erd], [BR], [Ri], [GP]).

\item A \textit{string module} St($w$) is associated to each string
representative $w$ by associating a quiver of type $\mathbb{A}_{n+1}$ to the
string $w$ with edges labeled by the letters $w_{i}$ and having the arrows
pointing to the left if $w$ is an arrow and to the right otherwise. The
\textquotedblleft obvious\textquotedblright\ representation of $Q$
corresponding to the diagram is specified as follows.

\item Let $w=\alpha_{n}^{\varepsilon_{n}}...\alpha_{1}^{\varepsilon_{1}}$ be a
string representative where the $\alpha_{i}$ are arrows and the $\varepsilon
_{i}$ are $-1$ or $1$. Consider the $\Bbbk$-vector space $V$ with basis
$v_{0},...,v_{n}.$ Define $V_{g}$ for each $g\in Q_{0}$ to be the span of the
following basis elements
\[
\{v_{i}|s(\alpha_{i+1}^{\varepsilon_{i+1}})=g\text{, }i=0,...,n-1\}
\]
together with $v_{n}$ if either: $t(\alpha_{n})=g$ and $\varepsilon_{n}=1,$ or
$s(\alpha_{n})=g$ and $\varepsilon_{n}=-1.$ In other words, we partition the
basis so that $v_{i}\in V_{s(\alpha_{i+1}^{\varepsilon_{i+1}})},$
$i=0,1,..n-1$ and $v_{n}\in V_{t(\alpha_{n}^{\varepsilon_{n}})}.$ The arrows
act by $\alpha_{i+1}(v_{i})=v_{i+1}$ if $\varepsilon_{i+1}=1$ and $\alpha
_{i}(v_{i})=v_{i-1}$ if $\varepsilon_{i}=-1$. All arrows not yet defined act
as zero. It is not hard to see that the representation does not depend on the
representative of the string, i.e. St$(w)\cong$St$(w^{-1})$. It is clear from
the construction that string modules are locally nilpotent and hence
correspond to $B$-comodules. We shall refer to such comodules as
\textit{string comodules}.

\item As an example consider the Kronecker quiver
\[
h\overset{a}{\underset{b}{\leftleftarrows}}g
\]
with arrows $a,b$ and vertices $g,h$ as shown. The string $b^{-1}ab^{-1}a$
corresponds the $\mathbb{A}_{5}$ quiver $v_{4}\rightarrow v_{3}\leftarrow
v_{2}\rightarrow v_{1}\leftarrow v_{0}$ or more suggestively,
\[%
\begin{tabular}
[c]{lllllllll}%
$v_{4}$ &  &  &  & $v_{2}$ &  &  &  & $v_{0}$\\
& $\overset{b}{\searrow}$ &  & $\overset{a}{\swarrow}$ &  & $\overset
{b}{\searrow}$ &  & $\overset{a}{\swarrow}$ & \\
&  & $v_{3}$ &  &  &  & $v_{1}$ &  &
\end{tabular}
\ \ \ \ \
\]
or even more briefly $%
%TCIMACRO{\QATOP{4\quad2\quad0}{3\quad1}}%
%BeginExpansion
\genfrac{}{}{0pt}{}{4\quad2\quad0}{3\quad1}%
%EndExpansion
$ just using subscripts and omitting the arrows that are recoverable with the
obvious convention. Here the sink vertices correspond to the socle $V_{g},$
and the source vertices correspond to the top $V_{h}$. The $5$-dimensional
string module $V$ decomposes as $V_{g}\oplus V_{h}$ as a $\Bbbk$-vector space,
where $V_{g}$ is spanned by $v_{4},v_{2},v_{0}$ and $V_{h}$ is spanned by
$v_{3},v_{1}$ with the indicated action of the arrows. Here the sink vertices
correspond to the socle $V_{g},$ and the source vertices correspond to the top
$V_{h}.$\bigskip
\end{itemize}

The basic result of Gelf'and-Ponomarev [GP] as reiterated by Ringel [Ri] is

\begin{theorem}
Let $A=\Bbbk(Q,I)$ be a string algebra. Then every finite dimensional
indecomposable module is isomorphic to a string module or a band module.
\end{theorem}

The following proposition shows that string (special biserial) algebras and
string (resp. special biserial) coalgebras are dual notions and that being a
string (resp. special biserial) coalgebra is local property. This will be used
in the specific example of quantum $SL(2)$ in the next section.

\begin{proposition}
Let $B$ be a pointed coalgebra.\newline(a) If $B$ is finite-dimensional, then
$B$ is a special biserial coalgebra if and only the dual algebra $DB$ is a
special biserial algebra.\newline(b) If $B$ is finite-dimensional, then $B$ is
a string coalgebra if and only if the dual algebra $DB$ is a string
algebra.\newline(c) Every subcoalgebra of a pointed special biserial coalgebra
is a special biserial coalgebra.\newline(d) Every finite-dimensional
subcoalgebra of a pointed special biserial coalgebra is contained in a
finite-dimensional subcoalgebra that is a string coalgebra.
\end{proposition}

\begin{proof}
Consider the finite-dimensional pointed coalgebra $B$ and view it as an
admissible subcoalgebra of the path algebra of its quiver $\Bbbk Q$. Set
$I=B^{\perp_{\Bbbk^{Q}}}\vartriangleleft\Bbbk^{Q}$. Then $I$ is an admissible
ideal of the path algebra $\Bbbk^{Q}.$ We have $B=C^{\perp_{\Bbbk^{Q}}\perp}$
by e.g. [Abe], Lemma 2.2.1, and $DB\cong\Bbbk^{Q}/I$. Parts (a) and (b) follow
from the easy observation that for a path $p\in\Bbbk Q$, $p^{\ast}\notin I$ if
and only if $p\vdash B.$

To prove (c), let $B$ be a subcoalgebra of a pointed special biserial
coalgebra which is an admissible subcoalgebra of the path coalgebra $\Bbbk Q.$
Then \textit{a fortiori} $B$ satisfies S1 and S2 in the definition, so by
Lemma 1.2.2 $B$ is special biserial.

Part (d) follows using the observation that one can enlarge a
finite-dimensional special biserial subcoalgebra $B$ of the string coalgebra
$B^{\prime}$ by adjoining all paths $p$ in $B^{\prime}$ such that $p\vdash B.$
This yields a finite-dimensional string coalgebra containing $B$.
\end{proof}

\begin{proposition}
Every finite dimensional indecomposable comodule over a string coalgebra is
either a string comodule or a band comodule.
\end{proposition}

\begin{proof}
The statement follows using duality from the Theorem and Proposition above.
\end{proof}

\section{The quantized coordinate Hopf algebra}

Assume the base field $\Bbbk$ is of characteristic zero. Henceforth, let
$C=\Bbbk_{\zeta}[SL(2)]$ be the $q$-analog of the coordinate Hopf algebra of
$SL(2)$, where $q$ is specialized to a root of unity $\zeta$ of odd order
$\ell$. We study the category \emph{M}$_{f}^{C}$ of finite-dimensional
comodules of $C$, which is equivalent to the category of finite-dimensional
modules over the Lusztig hyperalgebra $U_{\zeta}$ modulo the relation
$K^{\ell}-1$ (so only type 1 modules are considered), cf. [Lu], [APW]. A
fundamental fact is the existence of a nondegenerate Hopf pairing $U_{\zeta
}\otimes C\rightarrow\Bbbk\mathrm{.}$ which yields a Hopf algebra isomorphism
$C\rightarrow U_{\zeta}^{0}$ (see [CK2, section 5] for further references)
where $^{0}$ denotes the finitary dual [Mo]. Quantized coordinate algebras can
be defined by taking appropriate duals of quantized enveloping algebras, and
as in [APW, Appendix], it can be shown in type $\mathbb{A}$ that the resulting
Hopf algebra agrees with other definitions in the literature.

The coalgebra $C$ has the following presentation. The $\Bbbk$-algebra
generators are $a,b,c,d,$ with relations
\begin{align*}
ba  &  =\zeta ab\\
db  &  =\zeta db\\
ca  &  =\zeta ac\\
bc  &  =cb\\
ad-da  &  =(\zeta-\zeta^{-1})bc\\
ad-\zeta^{-1}bc  &  =1
\end{align*}
and with Hopf algebra structure further specified by
\begin{align*}
\Delta%
\begin{pmatrix}
a & b\\
c & d
\end{pmatrix}
&  =%
\begin{pmatrix}
a & b\\
c & d
\end{pmatrix}
\otimes%
\begin{pmatrix}
a & b\\
c & d
\end{pmatrix}
\\
\varepsilon%
\begin{pmatrix}
a & b\\
c & d
\end{pmatrix}
&  =%
\begin{pmatrix}
1 & 0\\
0 & 1
\end{pmatrix}
,\text{ }S%
\begin{pmatrix}
a & b\\
c & d
\end{pmatrix}
=%
\begin{pmatrix}
d & -\zeta b\\
-\zeta^{-1}c & a
\end{pmatrix}
\end{align*}

The injective indecomposable comodules and Gabriel quiver of $C$ were
determined in [CK] (see also [Ch]) and are summarized in the Theorem below.

For each nonnegative integer $r$, there is a unique simple module $L(r)$ of
highest weight $r$. These comodules exhaust the simple comodules ([APW], [Lu],
also [CP, 11.2]). We shall write
\[
r=r_{1}\ell+r_{0}%
\]
where $0\leq r_{0}<l$, the \textquotedblright short $p$-adic
decomposition\textquotedblright\ of $r$. Define an \textquotedblleft$\ell
$-reflection\textquotedblright\ $\tau:\mathbb{Z}\rightarrow\mathbb{Z}\ $by
\[
\tau(r)=(r_{1}+1)\ell+\ell-r_{0}-2
\]
if $r_{0}\neq\ell-1$, and $\tau(r)=r$ if $r_{0}=\ell-1$. Put $\sigma=\tau
^{-1}$ so that $\sigma(r)=(r_{1}-1)\ell+\ell-r_{0}-2$ for $r$ with $r_{0}%
\neq\ell-1$

Modular reduction from the generic Lusztig form yields the highest weight
module $V(r)$ with $\mathrm{\dim}V(r)=r+1$. When $r\in\mathbb{N}%
=\{0,1,2,...\}$ with $r_{0}\neq\ell-1$ and $r_{1}>0.$ this known as a
\textit{Weyl module} having socle series $%
%TCIMACRO{\QATOP{L(r)}{L(\sigma(r))}}%
%BeginExpansion
\genfrac{}{}{0pt}{}{L(r)}{L(\sigma(r))}%
%EndExpansion
.$ The \textit{dual Weyl module} has socle series $%
%TCIMACRO{\QATOP{L(\sigma(r))}{L(r)}}%
%BeginExpansion
\genfrac{}{}{0pt}{}{L(\sigma(r))}{L(r)}%
%EndExpansion
$for all such $r.$ On the other hand we have $V(r)=L(r)$ if $r=r_{0}$ or
$r_{0}=\ell-1.$ (cf. [CP, 11.2.7])

Lusztig's tensor product theorem [Lu] (see also [CP]) says that for all
$r\in\mathbb{N}$,
\[
L(r)\cong L(r_{0})\otimes L(r_{1}\ell)\cong L(r_{0})\otimes L_{U}%
(r_{1})^{\mathrm{Fr}}%
\]
where $L_{U}(r_{1})^{\mathrm{Fr}}$ is the twist by the quantum Frobenius map
\textrm{Fr} of the simple $U(sl(2))$ module $L_{U}(r_{1}).$

Set $I(r)=I(L(r))$ for all integers $r\geq0$, the injective envelope of the
simple comodule $L(r).$

\begin{theorem}
(a) ([APW]) If $r_{0}=\ell-1,\,$then $I(r)=L(r)$. \newline(b) ([CK]) If
$r<\ell-1$, then $I(r)$ has socle series with factors
\begin{gather*}
L(r)\\
L(\tau(r))\\
L(r)
\end{gather*}
\newline(c)([CK]) If $r\geq\ell$ (and $r_{0}\neq\ell-1$), then $I(r)$ has
socle series with factors
\begin{gather*}
L(r)\\
L(\sigma(r))\oplus L(\tau(r))\\
L(r)
\end{gather*}

\end{theorem}

\begin{proof}
(of (b) and (c)) We give a quick alternative to the approach in [CK].

As computed in [Ch] and [CK] the Cartan matrix for each nontrivial block is
the symmetric matrix $(c_{ij})$ indexed by $\mathbb{N}$ with nonzero entries
$c_{ii}=2$ and $c_{i,i+1}=1=$ $c_{i+1,i}$ for all $i\in\mathbb{N}.$ Fixing
$r=r_{0}<\ell-1,$ this is interpreted as saying that the multiplicity of
$L(\tau^{i}(r))$ in $I(\tau^{i+1}(r))$ equals $1$, the multiplicity of
$L(\tau^{i+1}(r))$ in $I(\tau^{i}(r))$ equals 1, and the multiplicity of
$L(\tau^{i}(r))$ in $I(\tau^{i}(r))$ equals $2$. Taking the socle series of
Weyl modules and dual Weyl modules into account, we see that the second socles
of the injective envelopes have simple factors as in the statements.

The desired conclusions will follow immediately once we show that there are no
nonsplit self-extensions of simples. This fact follows from the following
argument (cf. e.g. [Ja, p. 205]). Let $r\in\mathbb{N}.$ Suppose $0\rightarrow
L(r)\overset{j}{\rightarrow}E\overset{p}{\rightarrow}L(r)\rightarrow0$ is a
short exact sequence, and let $v\in L(r)$ be a weight vector of weight $r$
generating $L(r)$. Then there exists $v^{\prime}\in E$ of weight $r$ with
$p(v^{\prime})=v.$ Now let $N$ be the subcomodule of $E$ generated by
$v^{\prime}$. Since $N\nsubseteq\ker p,$ the simplicity of $i(L(r))$ forces
the sequence to split.
\end{proof}

\bigskip

This result determines the quiver as having vertices labelled by nonnegative
integers and with arrows $r\leftrightarrows s$ in case $r_{0}\neq l-1$ if and
only if $s=\tau(r)$ or $\tau(s)=r$. In case $r_{0}=l-1$, the simple block of
$r$ is a singleton (equivalently $L(r)$ is injective) and we call this block
``trivial''. Thus the nontrivial block containing $L(r),$ $r<l-1$ has quiver
\[
r\leftrightarrows\tau(r)\leftrightarrows\tau^{2}(r)\leftrightarrows\cdot
\cdot\cdot
\]
To simplify notation, we shall drop all but the exponents on $\tau$ and write
\[
0\leftrightarrows1\leftrightarrows2\leftrightarrows\cdot\cdot\cdot
\]

The injectives indecomposable comodules are all finite-dimensional, in
contrast to the injectives for the ordinary (nonquantum) modular coordinate
coalgebra. The result also shows that the coradical filtration is of length
2$.$

\bigskip

It follows that $C$ is a \textit{semiperfect} coalgebra ([Lin]) in the sense
that every finite dimensional comodule has a finite-dimensional injective
envelope (or, equivalently, projective cover) on both the left and the right.
Thus by [CKQ, Corollary 5.4], we know that the category $\mathcal{M}_{f}^{C}$
of finite-dimensional right comodules has almost split sequences.

Since $C$ is a Hopf algebra, we also have that $C$ is quasi-coFrobenius (or
"self-projective"), i.e. $C$ is projective as both a right and left
$C$-comodule. See [DNR] for discussion of these coalgebraic properties.

The next Lemma will be used in the proof of Proposition 2.7.1. Notation is as
in Lemma 2.05.

\begin{lemma}
$I(r)\cong I(r_{0})\otimes L(r_{1}\ell)$
\end{lemma}

\begin{proof}
We provide a direct proof by dimension count for the rank one case here. A
more general proof can be found in [APW2], 4.6.

We can write $r=r_{0}+r_{1}\ell$ with $r_{1}>0$ and $0\leq r_{0}<\ell-1.$ The
left hand side $I(r)$ has a Weyl module filtration $%
%TCIMACRO{\QATOP{V(r)}{V(\tau(r))}}%
%BeginExpansion
\genfrac{}{}{0pt}{}{V(r)}{V(\tau(r))}%
%EndExpansion
.$ Computing dimensions we have \textrm{dim}$I(r)$ $=r+1+(r_{1}+1)\ell
+\ell-r_{0}-2=2\ell(r_{1}+1).$ On the other hand, $I(r_{0})$ has filtration $%
%TCIMACRO{\QATOP{V(r_{0})}{V(\tau(r_{0}))}}%
%BeginExpansion
\genfrac{}{}{0pt}{}{V(r_{0})}{V(\tau(r_{0}))}%
%EndExpansion
$and \textrm{dim}$L(r_{1}\ell)=\dim L_{U}(r_{1})^{\mathrm{Fr}}=r_{1}+1$. Thus
the two sides have the same dimension. The Lusztig tensor product theorem
implies that $I(r_{0})\otimes L(r_{1}\ell)$ has $I(L(r))$ as a direct summand.
The desired conclusion follows.
\end{proof}

\subsection{The basic coalgebra}

Let $C_{r}$ denote a nontrivial block of $C$ as specified by a simple module
with highest weight $r=r_{0}<\ell-1$, i.e. $C_{r}$=$\bigoplus\limits_{i\in
\mathbb{N}}^{{}}I(\tau^{i}(r))^{\tau^{i}(r)+1}$ as above (exponent is the
multiplicity). By inspecting maps between indecomposable injectives we obtain

\begin{theorem}
The basic coalgebra $B$ of $C_{r}$ is the subcoalgebra of path coalgebra of
the quiver
\[
0\overset{b_{0}}{\underset{a_{0}}{\leftrightarrows}}1\overset{b_{1}}%
{\underset{a_{1}}{\leftrightarrows}}2\overset{b_{2}}{\underset{a_{2}%
}{\leftrightarrows}}\cdot\cdot\cdot
\]
spanned by the by group-likes $g_{i}$ corresponding to vertices and skew
primitive arrows $a_{i},b_{i}$ together with coradical degree two elements
\begin{align*}
d_{0}  &  :=b_{0}a_{0}\\
d_{i+1}  &  :=a_{i}b_{i}+b_{i+1}a_{i+1},\text{ }i\geq0.
\end{align*}

\end{theorem}

\begin{proof}
Let $E$ denote the basic injective $\bigoplus\limits_{i=0}^{\infty}I(\tau
^{i}(r))$ for $C_{r}$. Then the basic coalgebra $B\ $is defined to be the
coendomorphism coalgebra $h_{C}(E,E)$ and we see easily that $B=\bigoplus
\limits_{i=0}^{\infty}\mathrm{D}\operatorname*{Hom}(I(\tau^{i}(r)),E).$ Set
$I_{n}=I(\tau^{n}(r))$ for all $n\in\mathbb{N}$ with the convention that
$I_{n}=0$ if $n<0$. From the description of the indecomposable injectives
above, we see that
\[
\operatorname*{Hom}(I_{n},E)=\operatorname*{Hom}(I_{n},I_{n-1}\oplus
I_{n}\oplus I_{n+1}).
\]
By inspecting socle series we observe that $\operatorname*{Hom}(I_{n+1}%
,I_{n})$ and $\operatorname*{Hom}(I_{n},I_{n+1})$ are both one-dimensional.
Also $\operatorname*{Hom}(I_{n},I_{n})$ is two-dimensional, spanned by the
identity map and a map with kernel rad$(I_{n})$. Moreover it is easy to see
that there is a basis of $\operatorname*{Hom}(I_{n},E)$ consisting of maps
\begin{align*}
\alpha_{n}  &  \in\operatorname*{Hom}(I_{n},I_{n+1})\\
\beta_{n}  &  \in\operatorname*{Hom}(I_{n+1},I_{n})\\
\gamma_{n},\delta_{n}  &  \in\operatorname*{Hom}(I_{n},I_{n})
\end{align*}
such that
\begin{align*}
\gamma_{n}^{2}  &  =\gamma_{n}\\
\gamma_{n+1}\alpha_{n}  &  =\alpha_{n}=\alpha_{n}\gamma_{n},\\
\gamma_{n}\beta_{n}  &  =\beta_{n}=~\beta_{n}\gamma_{n+1}\\
\beta_{n}\alpha_{n}~  &  =\delta_{n}%
\end{align*}
for all $n\geq0.$ Also it is easy to see inductively that we may choose
$\beta_{n}$ (adjusting scalars) so that
\[
\delta_{n}=\alpha_{n-1}\beta_{n-1}%
\]
for all $n\geq1.$ All other products are zero. We choose dual basis elements
$\alpha_{n}^{\ast}\in\mathrm{D}\operatorname*{Hom}(I_{n},I_{n+1}),$ $\beta
_{n}^{\ast}\in\mathrm{D}\operatorname*{Hom}(I_{n+1},I_{n})$ and $\gamma
_{n}^{\ast}$, $\delta_{n}^{\ast}\in\mathrm{D}\operatorname*{Hom}(I_{n}%
,I_{n}),$ all in $\bigoplus\limits_{i,j}\mathrm{D}\operatorname*{Hom}%
(I_{i},I_{j})=B$.

It is straightforward to check that the elements $\alpha_{n},\beta_{n}%
,\gamma_{n}$ satisfy the following comultiplications.
\begin{align*}
\Delta\gamma_{n}^{\ast}  &  =\gamma_{n}^{\ast}\otimes\gamma_{n}^{\ast}.\\
\Delta\alpha_{n}^{\ast}  &  =\gamma_{n+1}^{\ast}\otimes\alpha_{n}^{\ast
}+\alpha_{n}^{\ast}\otimes\gamma_{n}^{\ast}\\
\Delta\beta_{n}^{\ast}  &  =\gamma_{n}^{\ast}\otimes\beta_{n}^{\ast}%
+\beta_{n+1}^{\ast}\otimes\gamma_{n}^{\ast}\\
\Delta\delta_{n+1}^{\ast}  &  =\gamma_{n+1}^{\ast}\otimes\delta_{n+1}^{\ast
}+\beta_{n+1}^{\ast}\otimes\alpha_{n+1}^{\ast}+\alpha_{n}^{\ast}\otimes
\beta_{n}^{\ast}+\delta_{n+1}^{\ast}\otimes\gamma_{n+1}^{\ast}\\
\Delta\delta_{0}^{\ast}  &  =\gamma_{0}^{\ast}\otimes\delta_{0}^{\ast}%
+\beta_{0}^{\ast}\otimes\alpha_{0}^{\ast}+\delta_{0}^{\ast}\otimes\gamma
_{0}^{\ast}.
\end{align*}
for all $n\geq0.$ For example, by evaluating $\delta_{n+1}^{\ast}$ on the
products $\beta_{n+1}\alpha_{n+1}=\alpha_{n}\beta_{n}=\delta_{n}$ and
$\delta_{n+1}\gamma_{n+1}=\delta_{n+1}=\delta_{n+1}\gamma_{n+1}$ we get the
expression for $\Delta\delta_{n+1}^{\ast}.$

An obvious embedding of $B$ into the path coalgebra is given as follows (cf.
[Ni]). Embed the degree one term by sending $\gamma_{n}^{\ast}$ to $g_{n},$
$\alpha_{n}^{\ast}$ to $a_{n}$, $\beta_{n}^{\ast}$ and (necessarily)
$\delta_{n}^{\ast}$ to $d_{n}.$
\end{proof}

\subsection{Indecomposables for the associated string coalgebra}

From the description of the nontrivial basic block $B=B_{r}$ above it is
immediate that $B$ is a special biserial coalgebra. The indecomposable
injective right $B$-comodules are
\begin{align*}
I_{n}  &  =\Bbbk g_{n}+\Bbbk a_{n-1}+\Bbbk b_{n}+\Bbbk d_{n}\\
I_{0}  &  =\Bbbk g_{0}+\Bbbk b_{0}+\Bbbk d_{0}%
\end{align*}
$(n\geq1)$ and they are also projective. Clearly rad$I_{n}=\Bbbk g_{n}+\Bbbk
a_{n-1}+\Bbbk b_{n}$ and rad$I_{0}=\Bbbk g_{0}+\Bbbk b_{0}$. So by deleting
all the basis elements $d_{n},$ $n\geq0$ we reduce to an associated string
coalgebra $B^{\prime}$ spanned by the arrows and group-likes $g_{i},$ and
arrows $a_{n},b_{n}$, as in 1.4.

By 1.4, the AR-quiver of the string coalgebra $B^{\prime}$ is the same as for
$B$ (and for $C_{r})$ except for the injective-projectives. The noninjective
indecomposable comodules for $B^{\prime}$ are determined by strings. The
simples correspond to the group-likes and the length two comodules (Weyl and
dual Weyl modules) are determined by arrows. More generally, indecomposables
are given by concatenations of arrows and their formal inverses. There are
strings of the form (up to equivalence)
\begin{align*}
&  a_{t}b_{t-1}^{-}...a_{s-1}b_{s+1}^{-}a_{s}\\
&  b_{t}^{-}a_{t-1}...b_{s-1}a_{s+1}b_{s}^{-}\\
&  a_{t}b_{t-1}^{-}...a_{s+1}b_{s}^{-}\\
&  b_{t}^{-}a_{t-1}...b_{s+1}^{-}a_{s}%
\end{align*}
where the subscripts form an interval $[s,t]$ of nonnegative integers strictly
increasing from right to left in the string . We adopt the following notation
for string modules
\begin{align*}
M(t,s)  &  =\mathrm{St}(a_{t-1}b_{t-2}^{-}...a_{s+1}b_{s}^{-})\\
M^{\prime}(t,s)  &  =\mathrm{St}(b_{t-1}^{-}a_{t-2}...b_{s+1}^{-}a_{s})\\
N(t,s)  &  =\mathrm{St}(a_{t-1}b_{t-2}^{-}...b_{s+1}^{-}a_{s})\\
N^{\prime}(t,s)  &  =\mathrm{St}(b_{t-1}^{-}a_{t-2}...a_{s+1}b_{s}^{-}),
\end{align*}
specifying each type of module by its starting and ending vertices $s$ and $t$
respectively. We declare that $M(i,i)=M^{\prime}(i,i)$ is the simple module
$S(i)=\mathrm{\Bbbk}g_{i}$ for all $i$.

\begin{proposition}
The comodules $M(t,s)$, $M^{\prime}(t,s)\,$with $t-s$ $\,$even, $\,$together
with the $N(t,s),N^{\prime}(t,s)$ with $t-s$ odd exhaust the indecomposable
$B^{\prime}$-comodules. The remaining $B$-comodules are the
projective-injectives $I_{n},$ $n\geq0.$
\end{proposition}

\begin{proof}
In view of 1.4 and 1.5 we just need to check that all strings for $B^{\prime}$
are of the form mentioned in the definition of the comodules in question. This
is simple verification.
\end{proof}

\begin{proposition}
The finite-dimensional indecomposable $B^{\prime}$-comodules can be
constructed as iterated pushouts and pullbacks of Weyl modules and dual Weyl modules.
\end{proposition}

\begin{proof}
The string comodule represented by $b_{t-1}^{-}$ (or $b_{t-1}$) has
composition length two and has structure described by the diagram
\[%
\begin{tabular}
[c]{lll}%
$t$ &  & \\
& $\searrow$ & \\
&  & $t-1$%
\end{tabular}
\ \ \ \ \ \ \ \
\]
where the Weyl module with highest weight $\sigma^{t}(r_{0})$ corresponds (via
Morita-Takeuchi equivalence) to $N^{\prime}(t,t-1)=\mathrm{St}(b_{t-1})$.
Dually, the string comodule represented by $a_{t}$ has structure described by
the diagram
\[%
\begin{tabular}
[c]{lll}
&  & $t$\\
& $\swarrow$ & \\
$t+1$ &  &
\end{tabular}
\ \ \ \ \ \ \ \
\]
where the dual Weyl comodule corresponds to $N(t+1,t)=\mathrm{St}(a_{t})$. The
string conmodule represented by $a_{t}b_{t-1}^{-}...a_{i-1}b_{s}^{-}$ is
constructed by the pullback
\[%
\begin{array}
[c]{ccc}%
\mathrm{St}(a_{t}b_{t-1}^{-}...a_{i-1}b_{s}^{-}) & \longrightarrow &
\mathrm{St}(a_{t})\\
\downarrow &  & \downarrow\\
\mathrm{St}(b_{t-1}^{-}...a_{i-1}b_{s}^{-}) & \longrightarrow & S(t)
\end{array}
\]
where we recall that $S(t)$ denotes the simple comodule at the vertex $t$. The
other string comodules can be constructed inductively similarly by pushouts or
pullbacks. Therefore, by category equivalence, the indecomposable $C$
-comodules can be thus constructed.
\end{proof}

\begin{corollary}
(a) The coalgebra $C=\Bbbk_{\zeta}[SL(2)]$ is tame of discrete comodule type
[Sim3].\newline(b) $C$ is the directed union of finite-dimensional coalgebras
of finite type.
\end{corollary}

\begin{proof}
(a) We must show that for each dimension vector $v\in K_{0}(C)$, there exist
only finitely many indecomposable right $C$-comodules $M$ with dimension
vector dim$M=v.$ If $M^{\prime\prime}$ is an indecomposable with
dim$M^{\prime\prime}=v,$ then $M^{\prime\prime}$ corresponds to one of the two
string comodules $M(t,s)$, $M^{\prime}(t,s)\,$\ (if $t-s$ $\ $is even)$\,$and
to one of $N(t,s),N^{\prime}(t,s)$ (if $t-s$ is odd), or $M^{\prime\prime}$ is
injective. This proves (a).

(b) Fix $n\geq1.$ Let $B_{n}$ denote the subcoalgebra of $B$ spanned by
\[
\{g_{i},a_{j},b_{j},d_{j}|0\leq i\leq n,0\leq j\leq n-1\}.
\]
$B_{n}$ is a subcoalgebra of path coalgebra of the finite quiver \newline%
\[
0\overset{b_{0}}{\underset{a_{0}}{\leftrightarrows}}1\overset{b_{1}}%
{\underset{a_{1}}{\leftrightarrows}}2\overset{b_{2}}{\underset{a_{2}%
}{\leftrightarrows}}\cdot\cdot\cdot n-1\overset{b_{n-1}}{\underset{a_{n-1}%
}{\leftrightarrows}}n.
\]
Clearly $B_{n}$ is a special biserial coalgebra, with dual special biserial
algebra $R_{n}=DB_{n}.$ The associated string coalgebra $B_{n}^{\prime}$ is
the span of $\{g_{i},a_{j},b_{j}|0\leq i\leq n,0\leq j\leq n-1\}.$ One sees
easily that the dual algebra $DB_{n}^{\prime}$ is a string algebra isomorphic
to the path algebra modulo the ideal of paths of length greater than one. Now
the strings are ones already determined for all of $B^{\prime},$ and, again,
there are no bands. So by the theory of special biserial algebras, see [Erd],
the string algebra $DB_{n}^{\prime}$ and the special biserial algebra $DB_{n}$
are both of finite representation type. \ By the standard category equivalence
between modules for a finite-dimensional algebra and comodules over its dual,
the same is true for the dual coalgebras $B_{n}^{\prime}$ and $B_{n}$.
Obviously $B=\cup B_{n}$, so the proof of (b) is finished.
\end{proof}

Part (b) is in contrast to the situation for restricted representations
([Xia], [Sut], [EGST], [RT]). Also $C$ is not locally of finite type in the
sense of Bongartz and Gabriel, see [BG].

\begin{proposition}
There is a duality on $\mathcal{M}_{f}^{B^{\prime}}$ which is obtained on
strings by interchanging the letters $a_{i}$ and $b_{i}^{-}$ for all i. The
duality interchanges $M(s,t)$ and $M^{\prime}(s,t),$ and interchanges $N(s,t)
$ and $N^{\prime}(s,t)$.
\end{proposition}

\begin{proof}
This is just the usual duality on $\mathcal{M}_{f}^{C}$ (using the antipode on
$C$), transported via category equivalence to $\mathcal{M}_{f}^{B}$ and
restricted to $\mathcal{M}_{f}^{B^{\prime}}.$ Alternatively one can observe
that interchanging the $a$'s and $b$'s yields a coalgebra anti-isomorphism of
$B$. The statements concerning the $M$'s and $N$'s are obvious from the definitions.
\end{proof}

\subsection{Syzygies and almost split sequences}

\begin{proposition}
Every finite-dimensional non-injective indecomposable $B$-comodule is in the
$\Omega^{\pm}$-orbit of some simple comodule.
\end{proposition}

\begin{proof}
This is verified by directly computing $\Omega^{-1}$:
\begin{align}
\Omega^{-1}(M(t,s))  &  =M(t+1,s-1)\text{ if }t>s>0\text{ (}t-s\text{
}even\text{) }\tag*{(1)}\label{1}\\
\Omega^{-1}(M(t,0))  &  =N(t+1,0),\text{ }t\geq0\tag*{(2)}\label{2}\\
\Omega^{-1}(N(t,s))  &  =N(t+1,s+1),\text{ }t>s\geq0\text{ (}t-s\text{ odd)}
\tag*{(3)}\label{3}%
\end{align}
This can be done as follows. The diagram
\[%
\begin{tabular}
[c]{lllll}
&  & $%
%TCIMACRO{\QATOP{\QTATOP{t}{t+1\text{ }t-1}}{t}}%
%BeginExpansion
\genfrac{}{}{0pt}{}{\genfrac{}{}{0pt}{1}{t}{t+1\text{ }t-1}}{t}%
%EndExpansion
\oplus%
%TCIMACRO{\QATOP{\QTATOP{t-2}{t-1\text{ }t-3}}{t-2}}%
%BeginExpansion
\genfrac{}{}{0pt}{}{\genfrac{}{}{0pt}{1}{t-2}{t-1\text{ }t-3}}{t-2}%
%EndExpansion
\oplus...\oplus%
%TCIMACRO{\QATOP{\QTATOP{s}{s+1\text{ }s-1}}{s}}%
%BeginExpansion
\genfrac{}{}{0pt}{}{\genfrac{}{}{0pt}{1}{s}{s+1\text{ }s-1}}{s}%
%EndExpansion
$ &  & \\
& $\overset{}{\nearrow}$ &  & $\overset{}{\searrow}$ & \\
$%
%TCIMACRO{\QATOP{t-1\quad...\quad s-2\quad s+1}{t\quad\quad t-2\quad...\quad
%s+2\quad s}}%
%BeginExpansion
\genfrac{}{}{0pt}{}{t-1\quad...\quad s-2\quad s+1}{t\quad\quad t-2\quad
...\quad s+2\quad s}%
%EndExpansion
$ &  & \multicolumn{1}{c}{$\overset{\Omega^{-1}}{\rightsquigarrow}$} &  & $%
%TCIMACRO{\QATOP{t\quad t-1\quad...\quad\quad\quad s}{t+1\quad t-2\quad
%\quad...\quad s+1\quad s-1}}%
%BeginExpansion
\genfrac{}{}{0pt}{}{t\quad t-1\quad...\quad\quad\quad s}{t+1\quad
t-2\quad\quad...\quad s+1\quad s-1}%
%EndExpansion
$%
\end{tabular}
\ \ \ \ \ \ \ \ \ \
\]
with the injective envelope on top and cokernel on the right$\ $demonstrates
$\Omega^{-1}(M(t,s))=M(t+1,s-1),t>s>1$ and similarly
\[%
\begin{tabular}
[c]{lllll}
&  & $%
%TCIMACRO{\QATOP{\QTATOP{t}{t+1\text{ }t-1}}{t}}%
%BeginExpansion
\genfrac{}{}{0pt}{}{\genfrac{}{}{0pt}{1}{t}{t+1\text{ }t-1}}{t}%
%EndExpansion
\oplus%
%TCIMACRO{\QATOP{\QTATOP{t-2}{t-1\text{ }t-3}}{t-2}}%
%BeginExpansion
\genfrac{}{}{0pt}{}{\genfrac{}{}{0pt}{1}{t-2}{t-1\text{ }t-3}}{t-2}%
%EndExpansion
\oplus...\oplus%
%TCIMACRO{\QATOP{\QTATOP{2}{3\text{ }1}}{2}}%
%BeginExpansion
\genfrac{}{}{0pt}{}{\genfrac{}{}{0pt}{1}{2}{3\text{ }1}}{2}%
%EndExpansion
\oplus%
%TCIMACRO{\QATOP{\QTATOP{0}{1}}{0}}%
%BeginExpansion
\genfrac{}{}{0pt}{}{\genfrac{}{}{0pt}{1}{0}{1}}{0}%
%EndExpansion
$ &  & \\
& $\overset{}{\nearrow}$ &  & $\overset{}{\searrow}$ & \\
$%
%TCIMACRO{\QATOP{t-1\quad...\quad\quad1}{t\quad t-2\quad...\quad2\quad0}}%
%BeginExpansion
\genfrac{}{}{0pt}{}{t-1\quad...\quad\quad1}{t\quad t-2\quad...\quad2\quad0}%
%EndExpansion
$ &  & \multicolumn{1}{c}{$\overset{\Omega^{-1}}{\rightsquigarrow}$} &  & $%
%TCIMACRO{\QATOP{\quad t\quad\quad...\quad\quad2\quad0}{t+1\quad t-1\quad
%...\quad\quad1\quad}}%
%BeginExpansion
\genfrac{}{}{0pt}{}{\quad t\quad\quad...\quad\quad2\quad0}{t+1\quad
t-1\quad...\quad\quad1\quad}%
%EndExpansion
$%
\end{tabular}
\ \ \ \ \ \ \ \ \ \
\]
yields $\Omega^{-1}(M(t,0))=N(t+1,0).$ Lastly, we obtain \ref{3}:
\[%
\begin{tabular}
[c]{lllll}
&  & $%
%TCIMACRO{\QATOP{\QTATOP{t}{t+1\text{ }t-1}}{t}}%
%BeginExpansion
\genfrac{}{}{0pt}{}{\genfrac{}{}{0pt}{1}{t}{t+1\text{ }t-1}}{t}%
%EndExpansion
\oplus%
%TCIMACRO{\QATOP{\QTATOP{t-2}{t-1\text{ }t-3}}{t-2}}%
%BeginExpansion
\genfrac{}{}{0pt}{}{\genfrac{}{}{0pt}{1}{t-2}{t-1\text{ }t-3}}{t-2}%
%EndExpansion
\oplus...\oplus%
%TCIMACRO{\QATOP{\QTATOP{s+1}{s+2\text{ \ }s}}{s+1}}%
%BeginExpansion
\genfrac{}{}{0pt}{}{\genfrac{}{}{0pt}{1}{s+1}{s+2\text{ \ }s}}{s+1}%
%EndExpansion
$ &  & \\
& $\overset{}{\nearrow}$ &  & $\overset{}{\searrow}$ & \\
$%
%TCIMACRO{\QATOP{\quad t-1\quad...\quad s+2\quad s}{t\quad t-2\quad...\quad
%s+1}}%
%BeginExpansion
\genfrac{}{}{0pt}{}{\quad t-1\quad...\quad s+2\quad s}{t\quad t-2\quad...\quad
s+1}%
%EndExpansion
$ &  & \multicolumn{1}{c}{$\overset{\Omega^{-1}}{\rightsquigarrow}$} &  & $%
%TCIMACRO{\QATOP{\quad t\quad...\quad\quad s+3\quad s+1}{t+1\quad
%t-1\quad...\quad s+2\quad}}%
%BeginExpansion
\genfrac{}{}{0pt}{}{\quad t\quad...\quad\quad s+3\quad s+1}{t+1\quad
t-1\quad...\quad s+2\quad}%
%EndExpansion
$%
\end{tabular}
\ \ \ \ \ \ \ \ \ \
\]

From the formulae above for $\Omega^{-1}$ we get the following expressions for
string comodules as elements in the $\Omega^{\pm}$-orbits of simples:%

\begin{align*}
M(t,s)  &  =\Omega^{-(\frac{t-s}{2})}(S(\frac{t+s}{2}))\\
N(t,s)  &  =\Omega^{-(\frac{t+s+1}{2})}(S(\frac{t-s-1}{2}))\\
M^{\prime}(t,s)  &  =\Omega^{\frac{t-s}{2}}(S(\frac{t+s}{2}))\\
N^{\prime}(t,s)  &  =\Omega^{\frac{t+s+1}{2}}(S(\frac{t-s-1}{2})).
\end{align*}
The first equation follows directly from \ref{1}. The second comes from
computing
\begin{align*}
&  \Omega^{-(\frac{t+s+1}{2})}(S(\frac{t-s-1}{2}))\\
&  =\Omega^{-(\frac{t+s+1}{2})}(M(\frac{t-s-1}{2},\frac{t-s-1}{2}))\\
&  =[\Omega^{-s}\circ\Omega^{-1}\circ\Omega^{-(\frac{t-s-1}{2})}%
](M(\frac{t-s-1}{2},\frac{t-s-1}{2}))\\
&  =[\Omega^{-s}\circ\Omega^{-1}](M(t-s-1,0))\\
&  =\Omega^{-s}(N(t-s,0))\\
&  =N(t,s),
\end{align*}
using \ref{1},\ref{2} and \ref{3} in succession. The second pair of equations
are obtained dually.
\end{proof}

The $\Omega^{\pm}$-orbits of indecomposable comodules for a nontrivial block
can be graphed by identifying the string comodules $M(x,y)$ and $N(x,y) $ with
the Cartesian points $(x,y)$ and identifying the dual comodules $M^{\prime
}(x,y)$ and $N^{\prime}(x,y)$ with the points $(y,x)$, reversing the
coordinates. The figure below shows the orbits containing the simples $(0,0)$
and $(3,3)$.%

%TCIMACRO{\FRAME{dpFX}{4.5254in}{3.0264in}{0pt}{}{}{jyzi440g.wmf}%
%width 4.5254in;  height 3.0264in;  depth 0pt;  original-width 4.583in;
%original-height 3.0556in;  cropleft "0";  croptop "1";  cropright "1";
%cropbottom "0";  filename 'JYZI440G.wmf';file-properties "XNPEU";}} }%
%BeginExpansion
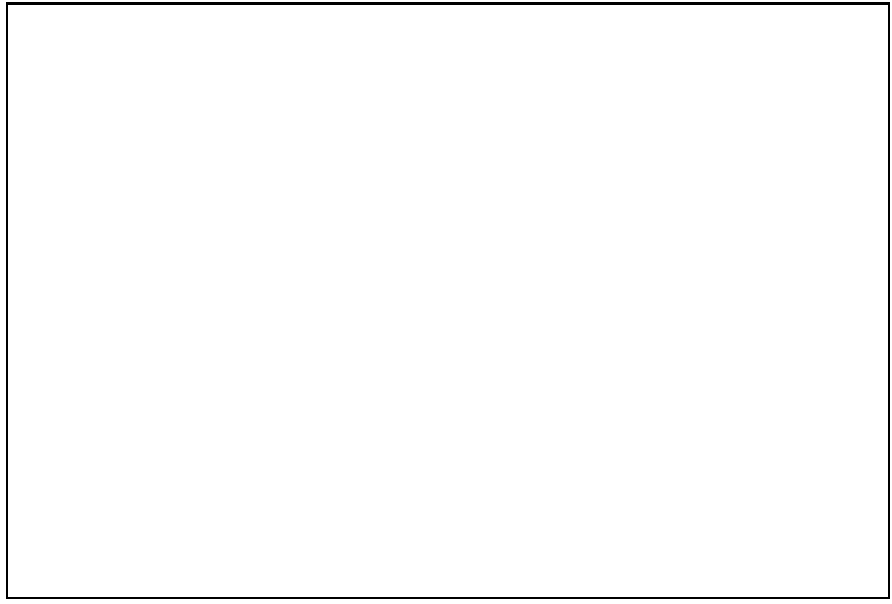
\begin{center}
\fbox{\includegraphics[
natheight=3.026400in,
natwidth=4.525400in,
height=3.0264in,
width=4.5254in
]%
{G:/graphics/JYZI440G__1.pdf}%
}
\end{center}
%EndExpansion

The almost split sequences and AR quiver for the category of
finite-dimensional $B$-comodules are described starting with the sequences
having an injective-projective comodule $I_{n}$ in the middle term. These are
precisely the sequences
\begin{equation}
0\rightarrow\mathrm{rad}(I_{n})\rightarrow\frac{\mathrm{rad}(I_{n}%
)}{\mathrm{soc}(I_{n})}\oplus I_{n}\rightarrow\frac{I_{n}}{\mathrm{soc}I_{n}%
}\rightarrow0 \tag*{(4)}\label{4}%
\end{equation}
with $n\in\mathbb{N}$. These sequences can be rewritten as
\begin{equation}
0\rightarrow\Omega(S(n))\rightarrow S(n-1)\oplus S(n+1)\oplus I_{n}%
\rightarrow\Omega^{-1}(S(n))\rightarrow0 \tag*{(5)}\label{5}%
\end{equation}
where $S(-1)$ is declared to be $0$.

\begin{theorem}
Applying $\Omega^{i}$, $i\in\mathbb{Z},$ to the sequences \ref{5} yields all
almost split sequences for $\mathcal{M}_{f}^{B}$
\end{theorem}

\begin{proof}
Since $B$ is self-projective (i.e. quasi-coFrobenius), the usual arguments for
modules (e.g. [ARS, IV.3]) apply to show that the functors $\Omega^{\pm1}$
provide autoequivalences of the stable category \underline{$\mathcal{M}$}%
$_{f}^{^{B}}$ (modulo injective-projectives). Applying $\Omega$ to the ends of
an almost split sequence $0\rightarrow M\rightarrow E\rightarrow
N\rightarrow0$ yields the exact commutative diagram
\[%
\begin{tabular}
[c]{ccccccccc}
&  & $0$ &  & $0$ &  & $0$ &  & \\
&  & $\downarrow$ &  & $\downarrow$ &  & $\downarrow$ &  & \\
$0$ & $\longrightarrow$ & $\Omega M$ & $\longrightarrow$ & $\Omega E\oplus J$
& $\longrightarrow$ & $\Omega N$ & $\longrightarrow$ & $0$\\
&  & $\downarrow$ &  & $\downarrow$ &  & $\downarrow$ &  & \\
$0$ & $\longrightarrow$ & $I_{M}$ & $\longrightarrow$ & $I_{M}\oplus I_{N}$ &
$\longrightarrow$ & $I_{N}$ & $\longrightarrow$ & $0$\\
&  & $\downarrow$ &  & $\downarrow$ &  & $\downarrow$ &  & \\
$0$ & $\longrightarrow$ & $M$ & $\longrightarrow$ & $E$ & $\longrightarrow$ &
$N$ & $\longrightarrow$ & $0$\\
&  & $\downarrow$ &  & $\downarrow$ &  & $\downarrow$ &  & \\
&  & $0$ &  & $0$ &  & $0$ &  &
\end{tabular}
\ \ \
\]
where $I_{M}\rightarrow M$ and $I_{N}\rightarrow N$ are projective covers and
$J$ is a projective-injective comodule. The arguments of [ARS, X.1,
Propositions 1.3, 1.4, p. 340] show that the top row $0\rightarrow\Omega
M\rightarrow\Omega E\oplus J\rightarrow\Omega N\rightarrow0$ is an almost
split sequence. Similarly, there is an almost split sequence $0\rightarrow
\Omega^{-1}M\rightarrow\Omega^{-1}E\oplus I\rightarrow\Omega^{-1}%
N\rightarrow0$ for some projective-injective $I$.

By the Proposition above, we obtain all almost split sequences in this manner.
Since the sequences in \ref{5} are precisely the almost split sequences with
an injective-projective in the middle term, the other sequences
\[
0\rightarrow\Omega^{i+1}(S(n))\rightarrow\Omega^{i}S(n-1)\oplus\Omega
^{i}S(n+1)\rightarrow\Omega^{i-1}(S(n))\rightarrow0
\]
we obtain when we apply $\Omega^{i},$ for a nonzero integer $i,$ do not have
an injective-projective summand in the middle term. Thus the sequences as in
the statement are all almost split and exhaust all almost split sequences.
\end{proof}

\subsection{The Auslander-Reiten quiver}

We define the Auslander-Reiten (AR) quiver of a coalgebra to be the quiver
whose vertices are isomorphism classes of indecomposable comodules and whose
(here multiplicity-free) arrows are defined by the existence of an irreducible
map between indecomposables. When the coalgebra $C$ is right semiperfect, then
the results of [CKQ] guarantee that the category of finite-dimensional right
$C$-comodules has almost split sequences, and thus that the AR quiver is
determined by them.

We now describe the AR quiver for $\mathcal{M}_{f}^{B}$.$\ $The stable AR
quiver (with all $I_{n}$'s deleted) consists of two components of type
$\mathbb{ZA}_{\infty}$ which are transposed by $\Omega$. To get the AR quiver
we use Proposition 1.2.5 (b). We adjoin the injective-projectives $I_{n}$ with
$n$ odd to the component containing $S(n)$ with $n$ odd, and similarly we
adjoin the $I_{n}$ with $n$ even to the other component containing the $S(n)$
with $n$ even. The two components for a nontrivial block are shown below.%

\begin{gather*}
\cdot\cdot\cdot%
\begin{array}
[c]{ccccccccccccccc}
& \Omega^{3}S(0) &  &  &  & \Omega S(0) & \rightarrow & I_{0} & \rightarrow &
\Omega^{-1}S(0) &  &  &  & \Omega^{-3}S(0) & \\
\nearrow &  & \searrow &  & \nearrow &  & \searrow &  & \nearrow &  & \searrow
&  & \nearrow &  & \searrow\\
&  &  & \Omega^{2}S(1) &  &  &  & S(1) &  &  &  & \Omega^{-2}S(1) &  &  & \\
\searrow &  & \nearrow &  & \searrow &  & \nearrow &  & \searrow &  & \nearrow
&  & \searrow &  & \nearrow\\
& \Omega^{3}S(2) &  &  &  & \Omega S(2) & \rightarrow & I_{2} & \rightarrow &
\Omega^{-1}S(2) &  &  &  & \Omega^{-3}S(2) & \\
\nearrow &  & \searrow &  & \nearrow &  & \searrow &  & \nearrow &  & \searrow
&  & \nearrow &  & \searrow\\
&  &  & \Omega^{2}S(3) &  &  &  & S(3) &  &  &  & \Omega^{-2}S(3) &  &  & \\
\searrow &  & \nearrow &  & \searrow &  & \nearrow &  & \searrow &  & \nearrow
&  & \searrow &  & \nearrow\\
& \Omega^{3}S(4) &  &  &  & \Omega S(4) & \rightarrow & I_{4} & \rightarrow &
\Omega^{-1}S(4) &  &  &  & \Omega^{-3}S(4) &
\end{array}
\cdot\cdot\cdot\\
\cdot\\
\cdot\\
\cdot
\end{gather*}

\begin{gather*}
\cdot\cdot\cdot%
\begin{array}
[c]{ccccccccccccccc}
& \Omega^{2}S(0) &  &  &  & S(0) &  &  &  & \Omega^{-2}S(0) &  &  &  &
\Omega^{-4}S(0) & \\
\nearrow &  & \searrow &  & \nearrow &  & \searrow &  & \nearrow &  & \searrow
&  & \nearrow &  & \searrow\\
&  &  & \Omega S(1) & \rightarrow & I_{1} & \rightarrow & \Omega^{-1}S(1) &  &
&  & \Omega^{-3}S(1) &  &  & \\
\searrow &  & \nearrow &  & \searrow &  & \nearrow &  & \searrow &  & \nearrow
&  & \searrow &  & \nearrow\\
& \Omega^{2}S(2) &  &  &  & S(2) &  &  &  & \Omega^{-2}S(2) &  &  &  &
\Omega^{-4}S(2) & \\
\nearrow &  & \searrow &  & \nearrow &  & \searrow &  & \nearrow &  & \searrow
&  & \nearrow &  & \searrow\\
&  &  & \Omega S(3) & \rightarrow & I_{3} & \rightarrow & \Omega^{-1}S(3) &  &
&  & \Omega^{-3}S(3) &  &  & \\
\searrow &  & \nearrow &  & \searrow &  & \nearrow &  & \searrow &  & \nearrow
&  & \searrow &  & \nearrow\\
& \Omega^{2}S(4) &  &  &  & S(4) &  &  &  & \Omega^{-2}S(4) &  &  &  &
\Omega^{-4}S(4) &
\end{array}
\cdot\cdot\cdot\\
\cdot\\
\cdot\\
\cdot
\end{gather*}

\begin{remark}
The almost split sequences with indecomposable middle terms are ones with
single arrows on the top boundary
\begin{align*}
0  &  \rightarrow\Omega^{i+1}(S(0))\rightarrow\Omega^{i}(S(1))\rightarrow
\Omega^{i-1}(S(0))\rightarrow0\text{ }\\
(0  &  \neq i\in\mathbb{Z)}%
\end{align*}
which are explicitly
\begin{align*}
0  &  \rightarrow N^{\prime}(i,i-1)\rightarrow M^{\prime}(i+2,i))\rightarrow
N^{\prime}(i+2,i+1)\rightarrow0\text{, }\\
0  &  \rightarrow N(i+2,i+1)\rightarrow M(i+2,i))\rightarrow
N(i,i-1)\rightarrow0\text{.}%
\end{align*}

\end{remark}

\subsection{Quantum dimension}

The quantum dimension of a $C$-comodule $M$ was defined in [And] to be the
trace of the action of $K$ as a linear transformation on $M$, i.e.,
\[
\operatorname*{qdim}M=\mathrm{Tr}_{K}(M)
\]

\begin{proposition}
Let $0\leq r_{0}<\ell-1$ and let $M_{r_{0}}(s,t)$ (resp. $N_{r_{0}}(s,t)$)
denote the $C$-comodule in the block of $L(r_{0})$ corresponding to the string
comodule $M(s,t)$ (resp. $N(s,t)$) The quantum dimensions are
\begin{align*}
\operatorname*{qdim}M_{r_{0}}(s,t)  &  =[r_{0}+1]_{\zeta}\sum_{i=s}%
^{t}(-1)^{i}(i+1)\\
\operatorname*{qdim}N_{r_{0}}(s,t)  &  =[r_{0}+1]_{\zeta}\sum_{i=s}%
^{t}(-1)^{i}(i+1).
\end{align*}

\end{proposition}

\begin{proof}
The composition series of $M_{r_{0}}(s,t)$ is $\{L(\sigma^{i}(r_{0}))|s\leq
i\leq t\}$ and is of odd length, The comodule $N_{r_{0}}(s,t)$ is of even
length has composition series given by the same expression. We compute
$\operatorname*{qdim}L(\sigma^{i}(r_{0}))=\mathrm{Tr}_{K}(L(\sigma^{i}%
(r_{0}))$ by first observing that $\sigma(r_{0})=\ell-r_{0}-2+\ell,$
$\sigma^{2}(r_{0})=r_{0}+2\ell$ and more generally
\[
\sigma^{j}(r_{0})=%
%TCIMACRO{\QATOPD{\{}{.}{\ell-r_{0}-2+j\ell\text{ if }j\text{ odd}}{r_{0}%
%+j\ell\text{ if }j\text{ even}}}%
%BeginExpansion
\genfrac{\{}{.}{0pt}{}{\ell-r_{0}-2+j\ell\text{ if }j\text{ odd}}{r_{0}%
+j\ell\text{ if }j\text{ even}}%
%EndExpansion
\text{.}%
\]

Note that $[\ell-r_{0}-2]_{\zeta}=-[r_{0}+1]_{\zeta}$. By the Lusztig tensor
product theorem (see [Lu] or [CP]) we have $L(r_{0}+j\ell)\cong L(r_{0}%
)\otimes L_{U}(j)^{\mathrm{Fr}}$ where $L_{U}(j)$ is the classical nonquantum
simple module of highest weight $j$ and the superscript $^{\mathrm{Fr}}$ is
the quantum Frobenius twist. The grouplike $K$ acts diagonally on the tensor
product and acts trivially on the second factor. Therefore
$\operatorname*{qdim}(L(r_{0}+j\ell))=(j+1)[r_{0}+1]_{\zeta}.$ Thus
\[
\operatorname*{qdim}L(\sigma^{j}(r_{0}))=%
%TCIMACRO{\QATOPD{\{}{.}{-(j+1)[r_{0}+1]_{\zeta}\text{ if }j\text{ odd}%
%}{(j+1)[r_{0}+1]_{\zeta}\text{ if }j\text{ even}}}%
%BeginExpansion
\genfrac{\{}{.}{0pt}{}{-(j+1)[r_{0}+1]_{\zeta}\text{ if }j\text{
odd}}{(j+1)[r_{0}+1]_{\zeta}\text{ if }j\text{ even}}%
%EndExpansion
\text{.}%
\]
The result follows.
\end{proof}

\begin{proposition}
The finite-dimensional noninjective indecomposable $C$-comodules are of
nonzero quantum dimension.
\end{proposition}

\begin{proof}
The result follows from the trivial observation that the alternating sum of a
nonempty sequence of consecutive integers is nonzero.
\end{proof}

\subsection{The stable Green ring}

The Green ring for $\mathcal{M}_{f}^{C}$ has as a $\mathbb{Z}$-basis of
isomorphism classes $[L]$ of finite dimensional indecomposable comodules
$L\in\mathcal{M}_{f}^{C}$, with addition given by direct sum and
multiplication given by the tensor product. The \textit{stable Green ring} is
the Green ring modulo the ideal of projective-injectives. This ring was
studied for certain finite-dimensional Hopf algebras in [EGST]. While the
basic coalgebra $B$ is not a Hopf algebra, $\mathcal{M}_{f}^{B}$ inherits a
quasi-tensor structure from $C$ (and from the quantized hyperalgebra
$U_{\zeta}.$).

By a well-known result [CP, 11.3] (cf. [And], see also [CK]) we have a formula
for the tensor product of two simples as a direct sum of simples and
injectives. It is used in the proof of the next result. We shall also use
\textit{Chebyshev polynomials of the second kind}. They are defined to be
polynomials $U_{n}(t)$ given recursively by
\[
U_{n}(t)=2tU_{n-1}(t)-U_{n-2}(t),
\]
with $U_{1}(t)=2t$, $U_{0}(t)=1.$ Note that $U_{n}(t)$ is of degree $n$ and is
an element of $\mathbb{Z}[2t]$.

\begin{proposition}
The stable Green ring is generated as a commutative $\mathbb{Z}$-algebra by
the isomorphism classes $x=[L(1)],$ $y=[L(\ell)],$ $w=[\Omega(L(0))]$ with
$w^{-1}=[$ $\Omega^{-1}(L(0))]$, subject to the relation $U_{\ell-1}(\frac
{1}{2}x).$
\end{proposition}

\begin{proof}
The Clebsch-Gordon decomposition for the tensor product of simple modules
$L(a)$ and $L(b)$ with for positive integers $a,b$ with $a\geq b$ is
\[
L(a)\otimes L(b)=%
%TCIMACRO{\QATOPD{\{}{.}{\dbigoplus _{t=a-b}^{a+b}L(t)\text{ if }a+b\leq\ell
%-1}{\dbigoplus \limits_{t=a-b}^{a^{\prime}+b^{\prime}}L(t)\oplus\dbigoplus
%\limits_{t=\ell-1}^{a+b}I(\sigma(t))\text{ if }a+b\geq\ell}}%
%BeginExpansion
\genfrac{\{}{.}{0pt}{}{{\displaystyle\bigoplus_{t=a-b}^{a+b}}L(t)\text{ if
}a+b\leq\ell-1}{{\displaystyle\bigoplus\limits_{t=a-b}^{a^{\prime}+b^{\prime}%
}}L(t)\oplus{\displaystyle\bigoplus\limits_{t=\ell-1}^{a+b}}I(\sigma(t))\text{
if }a+b\geq\ell}%
%EndExpansion
\]
where $t$ runs over integers such that $t\equiv a+b$ $\operatorname{mod}2,$
$a^{\prime}=\ell-a-2$ and $b^{\prime}=\ell-b-2.$ Let $x_{n}=[L(n)]$ and
$x=x_{1}$. This first case of the decomposition yields the recursive formula
\[
x_{n+1}=x\cdot x_{n}-x_{n-1}%
\]
for $n$ with $\ell-1>n>0$. Since $L(\ell-1)$ is projective, it is obvious that
$x_{\ell-1}=0.$ It follows that $x$ generates the subcoalgebra spanned by
$\{x_{i}|0\leq i\leq\ell-2\},$ and that $x$ satisfies the rescaled Chebyshev
polynomial of the second kind $f(x)=U_{\ell-1}(x/2)$ of degree $\ell-1.$

The simple modules $L(r)$ are of the form $L(r_{0})\otimes L(r_{1}\ell)\cong
L(r_{0})\otimes L_{U}(r_{1})^{\mathrm{Fr}}$ by the Lusztig tensor product
theorem [CP, 11.2], cf. [CK2]. Thus the remaining (nonrestricted) simples are
generated by $x$ along with the Frobenius twists of simple $U$-modules
$L_{U}(r_{1})^{\mathrm{Fr}}$. These modules are generated by $y:=L_{U}%
(1)^{\mathrm{Fr}}$, in view of the classical Clebsch-Gordon formula. Thus the
class of every simple module is uniquely of the form $x^{i}y^{j}$ where $0\leq
i\leq\ell-2$ and $j\geq0.$

Note that
\[
\Omega^{\pm}(M)\otimes N\cong M\otimes\Omega^{\pm}(N)\cong\Omega^{\pm
}(M\otimes N)
\]
for all $M,N\in\underline{\mathcal{M}}_{f}^{C}$. This says that in the stable
Green ring $\mathcal{R}$ the operator induced by $\Omega$, which we denote by
$\omega$, we have $\omega(xy)=x\omega(y)=\omega(x)y$ for all $x,y\in
\mathcal{R}$. Thus, by elementary ring theory, $\omega:\mathcal{R}%
\rightarrow\mathcal{R}$ is an $\mathcal{R}$-module map which equals
multiplication by
\[
w:=\omega(1_{\mathcal{R}})=[\Omega(L(0))]=[V(1)].
\]
Also $w$ is invertible with inverse $\omega^{-1}(1_{\mathcal{R}}%
)=[\mathrm{D}V(1)]$. For all $z\in R$ and integers $m,$ we now have
$\omega^{m}(z)=w^{m}z.$ Thus by Proposition 2.4.1 every indecomposable module
has image in $\mathcal{R}$ of the form $w^{m}z,$ $m\in\mathbb{Z}$ where $z$ is
the class of a simple module. By the proof of Proposition 2.4.1, all $\omega
$-orbits are infinite; thus we see that the class of every finite-dimensional
indecomposable comodule is uniquely represented by a monomial of the form
$w^{m}x^{i}y^{j}$ where $0\leq i\leq\ell-2$ and $j\geq0$.

It remains to see that the ideal $\mathcal{I}$ of injective-projectives of the
Green ring is generated by $[x_{\ell-1}]=U_{\ell-1}(\frac{1}{2}x).$ For the
non-obvious inclusion, let $I(L(r))$ be an indecomposable injective comodule
with socle $L(r)$. This injective can be written as $I(L(r_{0}))\otimes
L(r_{1}\ell)$ (Lemma 2.0.3). Thus it suffices to show that $I(L(r_{0}%
))\subset\mathcal{I}$. This follows from the decomposition
\[
L(\ell-1)\otimes L(b)=%
%TCIMACRO{\dbigoplus \limits_{t=\ell-1}^{a+b}}%
%BeginExpansion
{\displaystyle\bigoplus\limits_{t=\ell-1}^{a+b}}
%EndExpansion
I(\sigma(t))
\]
for $b=0,1,...,\ell-2$, which is special case of the formula at the beginning
of the proof. One observes iteratively that the classes of all summands which
occur on the right hand side (i.e. $I(\ell-1),$ $I(\ell)\oplus I(\ell-2),$
$I(\ell-1)\oplus I(\ell-3),...)$ are in $\mathcal{I}$.
\end{proof}

\begin{remark}
The conclusion that the Green ring of finite-dimensional modules is
commutative is immediate from the fact that M is a quasi-tensor category (cf.
[CP, p. 329]). But we do need to use this fact above.
\end{remark}

\end{document}